\documentclass[english]{IEEEtran4PSCC}
\usepackage[T1]{fontenc}
\usepackage[latin9]{inputenc}
\usepackage[cmex10]{amsmath}
\usepackage{amssymb}
\usepackage{graphicx,color}
\graphicspath{{./Figs/}}

\makeatletter

\providecommand{\tabularnewline}{\\}

\usepackage{algorithm, algorithmic}
\usepackage[noadjust]{cite}

\usepackage{caption}
\makeatother

\usepackage{babel}

\newcommand*{\TitleFont}{%
      \usefont{\encodingdefault}{\rmdefault}{n}{n}%
      \fontsize{16.5}{18.5}%
      \selectfont}

\IEEEoverridecommandlockouts

\begin{document}

\title{\TitleFont{Towards Optimal Energy Management of Microgrids with a Realistic Model}}

 \author{\IEEEauthorblockN{Wuhua Hu\IEEEauthorrefmark{1}$^{,\star}$\thanks{$^\star$ Since April 2016, W. Hu has been with the Signal Processing Department, Institute for
Infocomm Research, A*STAR, Singapore.},
 Ping Wang\IEEEauthorrefmark{2},
 Hoay Beng Gooi\IEEEauthorrefmark{1}}
 \IEEEauthorblockA{\IEEEauthorrefmark{1} School of Electrical and Electronic Engineering, Nanyang Technological University, Singapore}
 \IEEEauthorblockA{\IEEEauthorrefmark{2} School of Computer Engineering, Nanyang Technological University, Singapore}
\{hwh, wangping, ehbgooi\}@ntu.edu.sg
 }

\maketitle

\thispagestyle{empty}
\pagestyle{plain}

\begin{abstract}
This work considers energy management in a grid-connected microgrid which consists of multiple conventional generators (CGs), renewable generators (RGs)
and energy storage systems (ESSs). A two-stage optimization approach is presented to
schedule the power generation, aimed at minimizing the long-term average operating cost subject to operational and service constraints. The first stage of optimization
determines hourly unit commitment of the CGs via a day-ahead
scheduling, and the second stage performs economic dispatch of the CGs,
ESSs and energy trading via an hour-ahead scheduling.
The combined solution meets the need of handling large uncertainties in the load demand and renewable generation, and provides an efficient solution under limited computational resource which meets both short-term and long-term quality-of-service requirements. The performance of the
proposed strategy is evaluated by simulations based on real load demand and renewable generation data.
\end{abstract}


\section{Introduction}

Microgrids are considered as promising ingredients for smart grids
to integrate renewable energies in a reliable and economic manner.
As renewable resources are dynamic and stochastic in nature, the energy
management of microgrids becomes challenging and critical. This drives extensive research to resolve the challenges confronted \cite{liang2014stochastic}.
As a key step, it is important to establish a reliable model of the
microgrids that best describes the underlying physics. Various models
have been developed in the literature which often ignore or oversimplify certain complicate
operational constraints, such as the start-up/shut-down cost,
minimum on/off time and limited ramp up/down rate of conventional generators (CGs), the imperfect charging/discharging processes of energy storage systems (ESSs) and the associated aging cost, the forecast errors of the load demand and the generation of renewable generators (RGs). By using the simplified models, various
energy management systems have been developed in favor of their own
approaches while maintaining technical tractability \cite{parisio2011mixed,huang2013adaptive,lu2013online,zhang2013robust,sundistributed2015,shireal2015}.
However, this leaves questions on their actual performances because the effects of the
ignored complicate constraints are largely unexamined. Considering such a situation, it is of importance to re-think about the management system
design based on a more accurate microgrid model which takes
the realistic features and constraints into account.

This work makes an effort in this direction. We consider a
microgrid consisting of multiple RGs, CGs and ESSs, which
works in a grid-connected mode and can trade energy with the
external market. We develop
an energy management approach which first schedules the unit commitment
(UC) of the CGs in a day-ahead manner and then conducts economic dispatch
of the CGs, ESSs and energy trading in an hour-ahead
manner. The UC is aimed at minimizing the whole-day operating cost using
the day-ahead hourly forecast of the load demand, renewable supply and electricity
prices. The forecast errors could be large due to variations
and uncertainties in the demand and renewable energy generation. As a result, the solution obtained
may deviate much from the ideal optimal solution. This motivates
the real-time economic dispatch which maintains the UC solution but
uses more accurate one-hour forecast data to optimize the generation. Furthermore, the dispatch solution
approximately optimizes the long-term average operating cost of the
microgrid while meeting both a short-term and a long-term quality-of-service (QoS) requirement in addition
to the operational constraints. Because it only assumes unbiased and bounded
errors in the forecast data, the dispatch solution also avoids the need of identifying probability
distributions of the errors and enables a simple implementation.

Unlike the existing literature, the microgrid model we develop includes all typical components integrated in a flexible manner, and keeps almost all major features and realistic operational
constraints of each component, as briefly mentioned at the beginning. We also distinguish elastic from inelastic load demand, of which the former can be curtailed
or shifted while the latter must be satisfied. The long-term average curtailment is required to be below a given percentage. The mismatch between
the scheduled generation and the load demand is penalized, which
mimics the compensation for the degraded service to electricity users
and the adjustment cost induced to electricity providers. The joint consideration of these practical aspects leads to the new problem formulation and motivates the development of our unique two-stage approach for obtaining a solution.

The day-ahead UC and the hour-ahead economic dispatch
are formulated and solved as mixed-integer programs (MIPs) at the first and second stages, respectively. We set up a virtual microgrid and conduct simulations based on real trace data. The achieved economic performances are compared between the proposed strategy founded on a realistic microgrid model and its variants that assume
simplified models. The value of adopting a more realistic model is justified, and some discussions are made based on the new observations.

\textit{Notation.} We use $\mathcal{N}_{CG}$, $\mathcal{N}_{RG}$
and $\mathcal{N}_{ESS}$ to denote the sets of CGs, RGs and ESSs,
respectively. Given a scalar $x$, the term $(x)^{+}$ equals $\max\{x,0\}$. The time is discretized into slots as indicated by $0,1,...,T-1$ and collected in the set $\mathcal{T}$, where $T$ is the length of the planning
time horizon. Each time slot refers to a period of one hour. The words ``energy'' and
``power'' are used interchangeably whenever the energy is meant for a single
time slot. Boldface and normal letters are used to denote vectors
and scalars, respectively, and ``w.r.t.'' (or ``s.t.'') is an abbreviation for ``with respect
to'' (or ``subject to''). The rest notation is defined in the context when it first appears.

\section{Modeling the Microgrid}

This section models each component of a microgrid by describing its
cost function and operational constraints.

\subsection{CG Model}

Consider a CG $i\in\mathcal{N}_{CG}$ in time slot $t\in\mathcal{T}$.
Let its on/off status be indicated by 1/0 of $u_{i,t}$, and its power output
be scheduled as $p_{i,t}$. Its operating cost consists
of the start-up cost $c_{i}^{su}$, shut-down cost $c_{i}^{sd}$, fuel
cost $C_{i,t}^{f}(p_{i,t})$ and maintenance cost $C_{i,t}^{m}(p_{i,t})$,
of which the first two appear only if the relevant operations occur. The
total cost is thus expressed as
{\vspace{-12pt}\par\small{}
\begin{align*}
C_{i,t}(u_{i,t-1},u_{i,t},p_{i,t}) & =c_{i}^{su}u_{i,t}(1-u_{i,t-1})+c_{i}^{sd}u_{i,t-1}(1-u_{i,t})\\
 & \quad+C_{i,t}^{f}(p_{i,t})+C_{i,t}^{m}(p_{i,t}),
\end{align*}
}%
where we have assumed that $C_{i,t}^{m}(0)=0$ \cite{zhao2014short}.
The CG emits carbon gases when it is on, and the amount is estimated
as $E_{i}(p_{i,t})$. The actual amount is thus expressed as $u_{i,t}E_{i}(p_{i,t})$.

The operation of the CG must satisfy a couple of operational constraints.
Its output power is bounded as,
\begin{equation}
u_{i,t}p_{i,\min}\le p_{i,t}\le u_{i,t}p_{i,\max},\label{eq: CG power limits}
\end{equation}
where $p_{i,\min}$ and $p_{i,\max}$ are the lower and upper bounds, respectively. The power outputs of two sequential time slots cannot change too much, that is,
\begin{equation}
|p_{i,t}-p_{i,t-1}|\le r_{i}p_{i,\max}\label{eq: CG ramp rate}
\end{equation}
where $r_{i}\in(0,1]$ is the ramp rate coefficient.
Also, the CG must remain on (or off) in time slot $t$ if the duration
of the on (or off) status, denoted by $t_{i,t}^{\text{on}}$ (or $t_{i,t}^{\text{off}}$), is less than a given amount of time, denoted
by $T_{i,\min}^{\text{on}}$ (or $T_{i,\min}^{\text{off}}$). This calls for the following minimum on/off time constraints \cite{zhao2014short}:
\begin{equation}
\begin{aligned}\left(t_{i,t}^{\text{on}}-T_{i,\min}^{\text{on}}\right)\left(u_{i,t-1}-u_{i,t}\right) & \ge0,\\
\left(t_{i,t}^{\text{off}}-T_{i,\min}^{\text{off}}\right)\left(u_{i,t}-u_{i,t-1}\right) & \ge0.
\end{aligned}
\label{eq: CG minimum on-off times}
\end{equation}

\subsection{RG Model}

The RGs are assumed to be kept on and non-dispatchable. It is sufficient
to model their aggregate power output, denoted by $p_{RG,t}$ for
each time slot $t$, which is stochastic in time. Since the true value
of $p_{RG,t}$ for a future time slot $t$ is unknown, it is assumed to be estimated as $\hat{p}_{RG,t}$, which is unbiased
and has a bounded error as follows:
\begin{equation}
\mathbb{E}\{\hat{p}_{RG,t}\}=p_{RG,t},\thinspace\thinspace|\hat{p}_{RG,t}-p_{RG,t}|\le\delta_{RG,t},\label{eq: RG unbiased power estimate}
\end{equation}
where the symbol $\mathbb{E}$ is the expectation operator and $\delta_{RG,t}$ is a given bound of the forecast error.

\subsection{ESS Model}

Let the charge and discharge rates of ESS $i$ be scheduled
as $p_{i,t}^{c}$ and $p_{i,t}^{d}$ for time slot $t$, respectively, and let $v_{i,t}^{c}$ indicate the working mode of the ESS $i$, which is 1 (or 0) if it is \emph{not} discharged (or \emph{not} charged). These variables satisfy
\begin{equation}
0\le p_{i,t}^{c}\le v_{i,t}^{c} p_{i,\max}^{c},\thinspace\thinspace0\le p_{i,t}^{d}\le (1-v_{i,t}^{c}) p_{i,\max}^{d},\label{eq: charge-discharge incomptability and limits}
\end{equation}%
where $p_{i,\max}^{c}$ and $p_{i,\max}^{d}$ are the corresponding upper bounds. The two constraints ensure that charge and discharge comply with the rate limits and do not happen in the same time slot. After charge/discharge, the state
of charge (SOC) of the ESS, denoted by $s_{i,t}$, changes and renews
into
\begin{equation}
s_{i,t}=s_{i,t-1}+(\eta_{i}^{c}p_{i,t}^{c}-p_{i,t}^{d}/\eta_{i}^{d})/E_{i,\text{cap}},\label{eq: SOC update}
\end{equation}
where $\eta_{i}^{c},\eta_{i}^{d}\in(0,1]$ are the energy conversion
coefficients, and $E_{i,\text{cap}}$ is the energy capacity of the
ESS. To protect the ESS, the SOC must be maintained in the following
range:
\begin{equation}
s_{i,\min}\le s_{i,t}\le s_{i,\max},\label{eq: ES SOC limits}
\end{equation}
where $s_{i,\min},s_{i,\max}\in(0,1]$ are the given bounds.

Both charging and discharging processes reduce the ESS's life, incurring capital costs as estimated by $C_{i}^{c}(p_{i,t}^{c})$
and $C_{i}^{d}(p_{i,t}^{d})$, respectively, which are zero at zero
and increasing in the argument variables. The cost of operating
the $i$th ESS in time slot $t$ is thus given by
\[
C_{i}(p_{i,t}^{c},p_{i,t}^{d})=C_{i}^{c}(p_{i,t}^{c})+C_{i}^{d}(p_{i,t}^{d}).
\]

\subsection{Load Demand Model}

The demand is either inelastic or elastic load, which must be satisfied and can be curtailed/shifted, respectively. Let
their aggregate amounts be $d_{ie,t}$ and
$d_{e,t}$ in time slot $t$, respectively, satisfying $d_{e,t}\ge d_{e,\min}>0$ for a certain constant
$d_{e,\min}$. If part of the elastic demand is
unsatisfied, it incurs a shortage cost to compensate
the affected customers. Likewise, if the generation exceeds the total demand, it incurs a cost to adjust generation for balancing supply with demand.  The shortage (or surplus) cost is denoted by $C_{\text{short},t}(x_{t})$ (or $C_{\text{surp},t}(x_{t})$),
which is zero if the shortage (or surplus) amount $x_{t}$ is equal to 0.

Let the demands $d_{ie,t}$ and $d_{e,t}$ be forecasted as $\hat{d}_{ie,t}$
and $\hat{d}_{e,t}$, respectively, and the forecast errors satisfy the next conditions:
\begin{equation}
\begin{aligned}\mathbb{E}\{\hat{d}_{ie,t}\}=d_{ie,t}, & \thinspace\thinspace|\hat{d}_{ie,t}-d_{ie,t}|\le\delta_{ie,t},\\
\mathbb{E}\{\hat{d}_{e,t}\}=d_{e,t}, & \thinspace\thinspace|\hat{d}_{e,t}-d_{e,t}|\le\delta_{e,t},
\end{aligned}
\label{eq: Loads unbiased estimates}
\end{equation}
where $\delta_{ie,t}$ and $\delta_{e,t}$ are the given error bounds. For convenience, we also define the net demand as $d_{\text{net},t}\triangleq d_{ie,t}+d_{e,t}-p_{RG,t}$,
with a corresponding forecast value $\hat{d}_{\text{net},t}\triangleq\hat{d}_{ie,t}+\hat{d}_{e,t}-\hat{p}_{RG,t}$
and the associated error bound $\delta_{\text{net},t}\triangleq\delta_{ie,t}+\delta_{e,t}+\delta_{RG,t}$.
And let $\delta_{\text{net},\max}\triangleq\max_{t\in\mathcal{T}}\delta_{\text{net},t}$.

\subsection{Interaction with External Energy Market}

The microgrid is connected to a host power system and can trade energy
with the external energy market. Let the energy purchase and sale prices be given as $c_{b,t}$ and $c_{s,t}$
for time slot $t$, respectively, which satisfy $0<c_{p,\min}\le c_{p,t}\le c_{p,\max}$
and $0<c_{s,\min}\le c_{s,t}\le c_{s,\max}$. To avoid trading arbitrage,
we assume that $c_{p,t}>c_{s,t}$ (otherwise, the microgrid will always
purchase and then sell the same energy to make profit). Let the amount
of energy to purchase and sell be scheduled as $p_{p,t}$ and $p_{s,t}$,
respectively, as are bounded by,
\begin{equation}
0\le p_{p,t}\le p_{p,\max},\thinspace\thinspace0\le p_{s,t}\le p_{s,\max},\label{eq: buy-sell limits}
\end{equation}
where $p_{p,\max}$ and $p_{s,\max}$ are the transaction limits in
one time slot. The transaction incurs the following cost to the microgrid:
\[
C_{m,t}(p_{p,t},p_{s,t})=c_{p,t}p_{p,t}-c_{s,t}p_{s,t}.
\]
Minimization of this cost implies $p_{p,t}p_{s,t}=0$, i.e., purchase and sale of energy do not happen in the same time slot.

\subsection{Service Requirements}

The microgrid needs to meet some QoS requirements. The dispatchable
power generation from the CGs, ESSs and energy market
must meet the inelastic demand minus the non-dispatchable generation from the RGs, but
not exceed an upper bound of the net demand. Also, both the maximum
and the average percentages of unsatisfied elastic demand must
not exceed the levels given by $\alpha_{\max}$ and $\alpha_{\text{avg}}$, respectively.

Define the dispatchable supply as
\begin{equation}
p_{\text{disp},t}\triangleq\sum_{i\in\mathcal{N}_{CG}}p_{i,t}+\sum_{i\in\mathcal{N}_{ESS}}(p_{i,t}^{d}-p_{i,t}^{c})+p_{p,t}-p_{s,t}, \label{eq: disp_supply}
\end{equation}
and let $p_{t}\triangleq p_{\text{disp},t}+p_{RG,t}$, which is the
total energy generation. The aforementioned QoS requirements impose three constraints
for each time slot $t$: {$d_{ie,t}-p_{RG,t}\le p_{\text{disp},t}\le\hat{d}_{\text{net},t}+\delta_{\text{net},t}$,
$\frac{d_{ie,t}+d_{e,t}-p_{t}}{d_{e,t}}\le\alpha_{\max}$ and $\underset{T\rightarrow\infty}{\lim\sup}\frac{1}{T}\sum_{t=0}^{T-1}\mathbb{E}\left\{ \left(\frac{d_{ie,t}+d_{e,t}-p_{t}}{d_{e,t}}\right)^{+}\right\} \le\alpha_{\text{avg}}$},
where the two scalars satisfy $0\le\alpha_{\text{avg}}\le\alpha_{\max}<1$
(the time-average constraint becomes redundant if $\alpha_{\text{avg}}=\alpha_{\max}$).

By using forecast data, we reformulate the first two QoS constraints
and approximate the third one, yielding
{
\begin{equation}
\begin{aligned}p_{\text{disp},t} & \ge\underline{p}_{\text{disp},t}\triangleq\hat{d}_{\text{net},t}+\delta_{\text{net},t}-\alpha_{\max}(\hat{d}_{e,t}+\delta_{e,t}),\\
p_{\text{disp},t} & \le\bar{p}_{\text{disp},t}\triangleq\hat{d}_{\text{net},t}+\delta_{\text{net},t},
\end{aligned}
\label{eq: supply limit}
\end{equation}
\begin{equation}
\underset{T\rightarrow\infty}{\lim\sup}\frac{1}{T}\sum_{t=0}^{T-1}\mathbb{E}\left\{ \left(\frac{\hat{d}_{\text{net},t}-p_{\text{disp},t}}{\hat{d}_{e,t}}\right)^{+}\right\} \le\alpha_{\text{avg}},\label{eq: avg shortage limit}
\end{equation}
}%
where $\underline{p}_{\text{disp},t}$ and $\bar{p}_{\text{disp},t}$ are the lower and
upper bounds of the dispatchable generation, respectively, and in particular $\underline{p}_{\text{disp},t}$ is the tightest bound derived
by using (\ref{eq: RG unbiased power estimate}) and (\ref{eq: Loads unbiased estimates}). These constraints are more complicated than the QoS constraints used in the literature \cite{sundistributed2015,shireal2015}, because of taking the forecast errors of load demand and renewable generation into account.

In addition, the total carbon emission of the microgrid in each time
slot $t$ must not exceed a given amount $M_{t}$, namely,
\begin{equation}
\sum_{i\in\mathcal{N}_{CG}}u_{i,t}E_{i}(p_{i,t})\le M_{t},\label{eq: CG emission limits}
\end{equation}
and the microgrid must maintain a predefined operating reserve, $R_{t}$, to prepare for meeting emergent demand, i.e.,
\begin{equation}
\sum_{i\in\mathcal{N}_{CG}}(p_{i,\max}-p_{i,t})\ge R_{t}.\label{eq: CG operating reserves}
\end{equation}

\section{Problem Formulation}

The operation of the microgrid is reviewed periodically.
Each review aims to schedule the power generations of the CGs,
the charge/discharge rates of the ESSs and the trading with the external power market to minimize the time-average operating cost
of the microgrid, subject to constraints \eqref{eq: CG power limits}-\eqref{eq: CG operating reserves}. Some of the constraints are nonlinear, which make the optimization complex. Before presenting the optimization problem, we introduce linearizations of the nonlinear terms and constraints to enable a more efficient solution process.

\subsection{Linearization of Nonlinear Terms/Constraints}

First, introduce a binary variable $v_{i,t}$ satisfying $v_{i,t}=u_{i,t-1}u_{i,t}$,
for each $i\in\mathcal{N}_{CG}$ and $t\in\mathcal{T}$. Then the
CG cost {$C_{i,t}(u_{i,t-1},u_{i,t},p_{i,t})$} has an equivalent
form which is linear in $u_{i,t}$:
\begin{align*}
C_{i,t}(u_{i,t},v_{i,t},p_{i,t}) & =c_{i}^{su}(u_{i,t}-v_{i,t})+c_{i}^{sd}(u_{i,t-1}-v_{i,t})\\
 & \quad+C_{i}^{f}(p_{i,t})+C_{i}^{m}(p_{i,t}),
\end{align*}
where the symbol $C_{i,t}$ has been abused to denote the new
cost function. Under minimization of the above cost, the nonlinear
equality condition on $v_{i,t}$ is equivalent to
\begin{equation}
v_{i,t}\le u_{i,t-1},\,v_{i,t}\le u_{i,t},\label{eq: bilinear equivalence}
\end{equation}
which is a linear inequality for each $i\in\mathcal{N}_{CG}$ and
$t\in\mathcal{T}$.

Second, because $u_{i,t}-v_{i,t}$ and $u_{i,t-1}-v_{i,t}$ are indicators
of whether the CG $i$ is started up or shut down at the start of
time slot $t$, respectively, the minimum on/off time constraints
(\ref{eq: CG minimum on-off times}) can be reformulated as the
following linear inequalities:
\begin{equation}
\begin{aligned}u_{i,\tau} & \ge u_{i,t}-v_{i,t},\thinspace\thinspace\forall\,1\le\tau-t\le T_{i,\min}^{\text{on}}-1,\\
u_{i,\tau} & \le1-(u_{i,t-1}-v_{i,t}),\thinspace\thinspace\forall\,1\le\tau-t\le T_{i,\min}^{\text{off}}-1,
\end{aligned}
\label{eq: minimum on-off equivalence}
\end{equation}
for each $i\in\mathcal{N}_{CG}$ and $t\in\mathcal{T}$.

Third, the supply shortage and surplus cost functions can be written
more explicitly as $C_{\text{short},t}(w_{t})$ and $C_{\text{surp},t}(p_{\text{disp},t}-\hat{d}_{\text{net},t}+w_{t})$,
where $w_{t}$ is the amount of shortage satisfying:
\begin{equation}
w_{t}\ge 0,\,w_{t} \ge \hat{d}_{\text{net},t}-p_{\text{disp},t}.\label{eq: supply shortage}
\end{equation}
This is valid under the minimization of the two costs.

\subsection{The Energy Management Optimization Problem}

Let $\boldsymbol{u}_{t}\triangleq(u_{i,t})$ for all $i\in\mathcal{N}_{CG}$,
which is a column vector collecting the on/off status variables of
all CGs in time slot $t$, and $\boldsymbol{\pi}_{t}\triangleq[\boldsymbol{p}_{t}^{T},(\boldsymbol{p}_{t}^{c})^{T},(\boldsymbol{p}_{t}^{d})^{T},p_{p,t},p_{s,t}]^{T}$,
which collects all power variables of time slot $t$. The
vectors are defined as $\boldsymbol{p}_{t}\triangleq(p_{i,t})$ for
all $i\in\mathcal{N}_{CG}$, and $\boldsymbol{p}_{t}^{c}\triangleq(p_{i,t}^{c})$
and $\boldsymbol{p}_{t}^{d}\triangleq(p_{i,t}^{d})$ for all $i\in\mathcal{N}_{ESS}$.
We also let $\boldsymbol{z}_{t}\triangleq[w_{t},\boldsymbol{v}_{t}^{T},(\boldsymbol{v}_{t}^{c})^{T}]^{T}$
which collects all auxiliary variables, where $\boldsymbol{v}_{t}\triangleq(v_{i,t})$
for all $i\in\mathcal{N}_{CG}$, and $\boldsymbol{v}_{t}^{c}\triangleq(v_{i,t}^{c})$ for all $i\in\mathcal{N}_{ESS}$.
The energy management optimization problem is then defined as:
{\vspace{-12pt}\par\small{}
\begin{align*}
\textbf{P0:}\min_{\{\boldsymbol{\pi}_{t},\boldsymbol{z}_{t}\}_{t\in\mathcal{T}}} & \thinspace\thinspace\underset{T\rightarrow\infty}{\lim\sup}\dfrac{1}{T}\sum_{t=0}^{T-1}\mathbb{E}\left\{ J_{t}(\boldsymbol{u}_{t},\boldsymbol{\pi}_{t},\boldsymbol{z}_{t})\right\} \\
\mathrm{s.t.,} & \thinspace\thinspace\eqref{eq: CG power limits}-\eqref{eq: CG ramp rate},\eqref{eq: charge-discharge incomptability and limits}-\eqref{eq: ES SOC limits},\eqref{eq: buy-sell limits}-\eqref{eq: supply shortage}\\
 & \thinspace\thinspace u_{i,t},v_{i,t}\in\{0,1\},\thinspace\thinspace\forall \, i\in\mathcal{N}_{CG},t\in\mathcal{T},\\
 & \thinspace\thinspace v_{i,t}^{c}\in\{0,1\},\thinspace\thinspace\forall \, i\in\mathcal{N}_{ESS},t\in\mathcal{T},
\end{align*}
}%
where the operating cost is calculated by{\small{}
\begin{multline}
J_{t}(\boldsymbol{u}_{t},\boldsymbol{\pi}_{t},\boldsymbol{z}_{t})\triangleq\sum_{i\in\mathcal{N}_{CG}}\left[\begin{array}{c}
c_{i}^{su}(u_{i,t}-v_{i,t})+c_{i}^{sd}(u_{i,t-1}-v_{i,t})\\
+C_{i}^{f}(p_{i,t})+C_{i}^{m}(p_{i,t})
\end{array}\right]\\
+\sum_{i\in\mathcal{N}_{ESS}}\left[C_{i}^{c}(p_{i,t}^{c})+C_{i}^{d}(p_{i,t}^{d})\right]+c_{p,t}p_{p,t}-c_{s,t}p_{s,t}\\
+C_{\text{short},t}(w_{t})+C_{\text{surp},t}(p_{\text{disp},t}-\hat{d}_{\text{net},t}+w_{t}).\label{eq: total operation cost}
\end{multline}
}%
Problem \textbf{P0} has a time-average cost function, a time-average load
service constraint and a couple of operational constraints with real
and binary decision variables. These all make the problem hard to be solved exactly. A sub-optimal solution is thus explored in the next section which integrates an offline and an online solution.

\section{Approximate Solution Approach}

The start-up and shut-down costs and time coupling constraints make it difficult to optimize a CG's generation online. This motivates us to determine its hourly on/off status (i.e., the UC solution) day ahead, and to optimize the specific supplies from the CGs, ESSs and external market hour ahead while respecting the UC solution. This leads to a two-stage optimization approach for obtaining a sub-optimal solution to \textbf{P0}. The first-stage optimization is solved as an MIP with a 24-hour horizon, and the second-stage is reduced to an MIP with a one-hour horizon by applying the Lyapunov optimization method. Each of the two MIPs can be an MILP or MINLP, depending on the sub-cost functions in \eqref{eq: total operation cost}.

\subsection{Day-ahead UC}

The present purpose is
to determine the UC solution for each CG to be implemented in the
next day. With the forecast of electricity demand and renewable supply
available for the next day, we solve \textbf{P0} with $T=24$ h. The parameter $\alpha_{\max}$ is temporally set to the value of $\alpha_{\text{avg}}$ for computing a conservative schedule to tackle the uncertainties in the forecast data. In the optimization, we ensure that constraints
(\ref{eq: CG ramp rate}), (\ref{eq: SOC update}), (\ref{eq: bilinear equivalence})
and (\ref{eq: minimum on-off equivalence}), which interlace variables
of two neighboring days, are correctly enforced by using the values
of relevant variables realized in the last day. The resulting problem is an MIP
and can be solved offline by standard solvers.

Suppose the full solution is obtained as $\{\boldsymbol{u}_{t}^{*},\boldsymbol{\pi}_{t}^{*},\boldsymbol{z}_{t}^{*}\}_{t\in\mathcal{T}}$
for a 24 h of interest. We then keep the UC solution
$\{\boldsymbol{u}_{t}^{*}\}_{t\in\mathcal{T}}$, but re-optimize the
remaining solution by using more accurate one-hour forecast data in the next stage.

\subsection{Real-time Economic Dispatch}

Because of variations and uncertainties in the demand and renewable supply,
the energy generation scheduled in the first stage may deviate much
from the optimal solution. To mitigate the deviation, we proceed to
re-optimize the energy generation using the more accurate hour-ahead forecast data while maintaining
the UC solution obtained in the first stage.

To respect
the UC solution obtained in the first stage, we need to enforce the following constraint in this real-time stage:
\begin{equation}
p_{i,t}\le\tau_{i,t}r_{i}p_{i,\max}, \label{eq: dispatch_feasiblity}
\end{equation}
if a CG $i$ has been scheduled to be on at time $t$ until $t+\tau_{i,t}$, where $\tau_{i,t}$ is the duration of the on status. This constraint ensures that it is feasible in the online stage for the CG to decrease its generation to zero at a time scheduled in the first stage, without violating the ramp rate constraint given in \eqref{eq: CG ramp rate}.

\subsubsection{Solution via Lyapunov optimization}

The Lyapunov optimization method is mainly developed to solve optimization
problems with a time-averaged objective and time-average equality/inequality
constraints, in which the decision variables are constrained in some
implicit sets that depend only on the present time slot \cite{neely2010stochastic}.
To apply the method, we ignore the ramp rate constraint (\ref{eq: CG ramp rate})
that crosses two time slots but will restore it later to enforce feasibility, and also
to relax the SOC constraint (\ref{eq: ES SOC limits}) as $\lim_{T\rightarrow\infty}\frac{1}{T}\mathbb{E}\{s_{i,T}-\beta_{i}\}=0$,
where $\beta_{i}>0$ and its value needs to be selected appropriately.

Virtual queues for the time-averaged terms are carefully defined based on
the Lyapunov optimization method. To meet the time-average service
constraint (\ref{eq: avg shortage limit}), a virtual queue backlog
$Q_{t}$ is introduced which evolves as follows:
{\vspace{-12pt}\par\small{}
\begin{equation}
Q_{t}=\max\left\{ Q_{t-1}+\left(\frac{\hat{d}_{\text{net},t}-p_{\text{disp},t}}{\hat{d}_{e,t}}\right)^{+}-\alpha_{\text{avg}},0\right\} .\label{eq: Q iteration}
\end{equation}
}It accumulates the percentage of unsatisfied elastic demand.

Meanwhile, to meet the relaxed SOC constraint, which is a time-average
constraint on the charge and discharge rates, another virtual queue
backlog $S_{i,t}$ is defined for each $i\in\mathcal{N}_{ESS}$ which
evolves as follows:
\begin{equation}
S_{i,t}=S_{i,t-1}+q_{i,t}=s_{i,t}-\beta_{i},\label{eq: S definition}
\end{equation}
where $q_{i,t}\triangleq(\eta_{i}^{c}p_{i,t}^{c}-p_{i,t}^{d}/\eta_{i}^{d})/E_{i,\text{cap}}$ refers to the change in the SOC, and the parameter $\beta_{i}$
is selected such that the SOC value $s_{i,t}$ can be controlled within
the given limits without imposing the original SOC constraint.

Let
$c_{i,\max}^{c}\triangleq\max_{x\in[0,p_{i,\max}^{c}]}\frac{dC_{i}^{c}(x)}{dx}$
and $c_{i,\max}^{d}\triangleq\max_{x\in[0,p_{i,\max}^{d}]}\frac{dC_{i}^{d}(x)}{dx}$.
It can be shown that it is sufficient to select the parameter $\beta_{i}$
as
\begin{equation}
\beta_{i}=s_{i,\min}+\frac{p_{i,\max}^{d}}{\eta_{i}^{d}E_{i,\text{cap}}}+\frac{VE_{i,\text{cap}}}{\eta_{i}^{c}}\left(c_{i,\max}^{c}+c_{p,\max}\right)\label{eq: beta}
\end{equation}
for each $i\in\mathcal{N}_{ESS}$, where $0<V\le V_{\max}$ and $V_{\max}$
is given by
{\vspace{-12pt}\par\footnotesize{}
\begin{equation}
V_{\max}=\min_{i\in\mathcal{N}_{ESS}}\left\{ \frac{s_{i,\max}-s_{i,\min}-(\eta_{i}^{c}p_{i,\max}^{c}+p_{i,\max}^{d}/\eta_{i}^{d})/E_{i,\text{cap}}}{E_{i,\text{cap}}[(c_{i,\max}^{c}+c_{p,\max})/\eta_{i}^{c}+\eta_{i}^{d}(c_{i,\max}^{d}-c_{s,\min})]}\right\} .\label{eq: V_max}
\end{equation}
}The derivations of $\beta_{i}$ and $V_{\max}$ are motivated by
the recent work \cite{sundistributed2015}. The differences and
complications come from the imperfect charging and discharging processes
and the SOC expressed in percentage of the energy storage capacity.

Stack the two types of virtual queues into a vector, yielding $\Theta_{t}\triangleq[Q_{t},S_{i_{1},t},S_{i_{2},t},\dots,S_{i_{|\mathcal{N}_{ESS}|},t}]^{T}$
with $i_{1},i_{2},\dots,i_{|\mathcal{N}_{ESS}|}\in\mathcal{N}_{ESS}$.
Define the Lyapunov function as
{\vspace{0pt}\par\small
\begin{equation*}
L(\Theta_{t})  =  \frac{1}{2}\left(Q_{t}^{2}+\sum_{i\in\mathcal{N}_{ESS}}S_{i,t}^{2}\right),
\end{equation*}
}%
which is positive definite, since $L(\Theta_{t})>0$ if $\Theta_{t}\ne\boldsymbol{0}$
and $L(\Theta_{t})=0$ if and only if $\Theta_{t}=\boldsymbol{0}$.
The one-slot conditional Lyapunov drift is then obtained as
{\vspace{-9pt}\par\small{}
\begin{align*}
\Delta L(\Theta_{t}) & =\mathbb{E}\left\{ L(\Theta_{t})-L(\Theta_{t-1})\mid\Theta_{t-1}\right\} \\
 & =\frac{1}{2}\mathbb{E}\left\{ \left.Q_{t}^{2}-Q_{t-1}^{2}+\sum_{i\in\mathcal{N}_{ESS}}(S_{i,t}^{2}-S_{i,t-1}^{2})\right|\Theta_{t-1}\right\} ,
\end{align*}
}%
which can be shown to be bounded as
{\vspace{-6pt}\par\scriptsize{}
\begin{align}
\Delta L(\Theta_{t}) & \le B+\mathbb{E}\left\{ \left.\begin{array}{c}
Q_{t-1}\left(\left(\frac{\hat{d}_{\text{net},t}-p_{\text{disp},t}}{\hat{d}_{e,t}}\right)^{+}-\alpha_{\text{avg}}\right)\\
+\sum_{i\in\mathcal{N}_{ESS}}q_{i,t}S_{i,t-1}
\end{array}\right|\Theta_{t-1}\right\} ,\label{eq: drift bound}
\end{align}
}%
where the constant $B$ is given by
{\tiny{}
\begin{align*}
B & =\frac{1}{2}\max\left\{ \alpha_{\max}^{2},\frac{\delta_{\text{net},\max}^{2}}{d_{e,\min}^{2}}\right\} +\frac{1}{2}\sum_{i\in\mathcal{N}_{ESS}}\frac{\max\left\{ \left(\eta_{i}^{c}p_{i,\max}^{c}\right)^{2},\left(p_{i,\max}^{d}/\eta_{i}^{d}\right)^{2}\right\} }{E_{i,\text{cap}}^{2}}.
\end{align*}
}

Next, the following drift-plus-penalty cost is introduced to reflect
on the trade-off between constraint satisfaction and cost minimization,
whose upper bound is derived as follows:
{\vspace{-12pt}\par\small{}
\begin{multline}
\Delta L(\Theta_{t})+V\mathbb{E}\left\{ J_{t}(\boldsymbol{u}_{t}^{*},\boldsymbol{\pi}_{t},\boldsymbol{v}_{t})\mid\Theta_{t-1}\right\} \\
\le\text{right-hand-side of }\eqref{eq: drift bound}+V\mathbb{E}\left\{ J_{t}(\boldsymbol{u}_{t}^{*},\boldsymbol{\pi}_{t},\boldsymbol{v}_{t})\mid\Theta_{t-1}\right\} ,\label{eq: drift plus penalty bound}
\end{multline}
}%
where $V$ is the penalty scalar that controls the trade-off. The
Lyapunov optimization method minimizes the above upper bound
by greedily solving the following optimization problem for each time slot:
\begin{align*}
\textbf{P1:}\min_{\boldsymbol{\pi}_{t},\boldsymbol{z}_{t}} & \thinspace\thinspace VJ_{t}(\boldsymbol{u}_{t}^{*},\boldsymbol{\pi}_{t},\boldsymbol{z}_{t})+\sum_{i\in\mathcal{N}_{ESS}}q_{i,t}S_{i,t-1}+\frac{Q_{t-1}}{\hat{d}_{e,t}}w_{t}\\
\text{s.t.,} & \thinspace\thinspace\eqref{eq: CG power limits}-\eqref{eq: CG ramp rate},\eqref{eq: charge-discharge incomptability and limits}-\eqref{eq: ES SOC limits},\eqref{eq: buy-sell limits}-\eqref{eq: CG operating reserves},\eqref{eq: supply shortage},\eqref{eq: dispatch_feasiblity}\\
 & \thinspace\thinspace v_{i,t}^{c}\in\{0,1\},\thinspace\thinspace\forall \, i\in\mathcal{N}_{ESS},
\end{align*}
where the indicators $\{u_{i,t}\}$ in the constraints are fixed to
the corresponding values of $\{u_{i,t}^{*}\}$ as obtained in the
first-stage optimization. The ramp rate constraint (\ref{eq: CG ramp rate})
 has been restored and the SOC constraint \eqref{eq: ES SOC limits} is kept to ensure a feasible solution even if the parameter $V$ is set to a value larger than $V_{\max}$.

Problem \textbf{P1} is an MIP with the binary variables only indicating charge/discharge
operations of the ESSs in a single time slot. The problem can be solved
efficiently using mature MIP solvers. Let the solution be obtained as $(\boldsymbol{\pi}_{t}^{\star},\boldsymbol{z}_{t}^{\star})$.

The performance gap to the optimal schedule can be explicitly
derived under the strong assumption that the inelastic and elastic demands and the RGs' aggregate generation are respectively identically and independently
distributed over all time slots, as was similarly done with a simplified microgrid model in \cite{sundistributed2015}.
Since in reality this assumption is not satisfied in general, we rely on simulations instead to evaluate the performance of the obtained
sub-optimal schedule.

\section{Performance Evaluation}

This section sets up a virtual microgrid and conducts simulations with real trace data. The objectives are two folds: (i) corroborating the empirical performance of the proposed energy management approach under
demand and supply uncertainties; and (ii) investigating the impacts
of employing oversimplified component models on the generation schedule and the economic performance.

\subsection{The Microgrid Setup}

The simulated microgrid has one wind farm, three CGs and two
ESSs, and it interacts with an external energy market. The specifications of these components are given below. The units of mass, time, power, energy and cost are in
kg, h, kW, kWh and \$ (US), respectively, unless otherwise specified.
The parameter $V$ is set to $V_{\max}$ as computed per (\ref{eq: V_max}).

\textbf{RGs: }The RGs mimic a group of wind turbines that form a wind
farm. The aggregate power trace is extracted from a wind farm in California
for one summer week \cite{windtrace2004},
and then scaled to a maximum output of 1.2 MW, as shown in Fig. \ref{fig: One-stage-vs-two-stage}.
The aggregate output power in each time slot is assumed to be forecasted
with a zero-mean error that follows a uniform distribution. The maximum error is proportional to the absolute change in the
output power of the time slot w.r.t. that of the previous time slot, and the proportion coefficient is increasing in the forecast horizon. Further, the error is capped within 80\% of the minimum and 120\% of the maximum output powers.

\textbf{CGs: }There are three CGs, with aggregate output power up
to 3 MW. Their parameters are given in Table \ref{tb: CG parameters} \cite{rahbar2015real}, and their cost and
carbon emission functions are given in Table \ref{tb: CG cost functions} \cite{mitigation2011ipcc}.

\begin{table}
\protect\caption{\textsc{\small{}Parameters of CGs}. \label{tb: CG parameters}}

\noindent \centering{}{\scriptsize{}}%
\begin{tabular}{cccccc}
\hline
{\scriptsize{}Type} & {\scriptsize{}$p_{i,\min}$} & {\scriptsize{}$p_{i,\max}$} & {\scriptsize{}$r_{i}$} & {\scriptsize{}$T_{i,\min}^{\text{on}}$} & {\scriptsize{}$T_{i,\min}^{\text{off}}$}\tabularnewline
\hline
{\scriptsize{}1} & {\scriptsize{}90} & {\scriptsize{}600} & {\scriptsize{}0.60} & {\scriptsize{}2} & {\scriptsize{}2}\tabularnewline
{\scriptsize{}2} & {\scriptsize{}200} & {\scriptsize{}1000} & {\scriptsize{}0.55} & {\scriptsize{}3} & {\scriptsize{}3}\tabularnewline
{\scriptsize{}3} & {\scriptsize{}350} & {\scriptsize{}1400} & {\scriptsize{}0.50} & {\scriptsize{}4} & {\scriptsize{}4}\tabularnewline
\hline
\end{tabular}
\end{table}

\begin{table}
\protect\caption{\textsc{\small{}Costs and carbon emissions of CGs}.\label{tb: CG cost functions}}

\noindent \centering{}{\scriptsize{}}%
\begin{tabular*}{1\columnwidth}{@{\extracolsep{\fill}}ccccc}
\hline
{\scriptsize{}Type} & {\scriptsize{}$c_{i}^{su}=c_{i}^{sd}$} & {\scriptsize{}$C_{i}^{f}(p_{i,t})$} & {\scriptsize{}$C_{i}^{m}(p_{i,t})$} & {\scriptsize{}$E_{i}(p_{i,t})$}\tabularnewline
\hline
{\scriptsize{}1} & {\scriptsize{}49.2} & {\scriptsize{}$1.72\times10^{-6}p_{i,t}^{2}+0.055p_{i,t}$} & {\scriptsize{}$0.026p_{i,t}$} & {\scriptsize{}$0.475p_{i,t}$}\tabularnewline
{\scriptsize{}2} & {\scriptsize{}79.7} & {\scriptsize{}$1.66\times10^{-6}p_{i,t}^{2}+0.053p_{i,t}$} & {\scriptsize{}$0.025p_{i,t}$} & {\scriptsize{}$0.472p_{i,t}$}\tabularnewline
{\scriptsize{}3} & {\scriptsize{}108.1} & {\scriptsize{}$1.59\times10^{-6}p_{i,t}^{2}+0.051p_{i,t}$} & {\scriptsize{}$0.024p_{i,t}$} & {\scriptsize{}$0.465p_{i,t}$}\tabularnewline
\hline
\end{tabular*}
\end{table}

\textbf{ESSs: }The two ESSs are composed of small Li-ion batteries, and their specifications are given in Table \ref{tb: ES parameters}.
The charging and discharging cost functions are approximated by three
piece-wise linear functions, which were experimentally obtained for
a Li-ion battery with a capacity of 8.1 W \cite{trippe2014charging}.
Each cost function accounts for the battery's life loss in terms of
the initial purchase cost. When the
operation is updated at a period of one hour, the cost has
the following explicit form:
{\par \vspace{-12pt} \small
\[
\begin{aligned}C_{i}(p_{i,t}^{c},p_{i,t}^{d})\approx \min_{z_i}\frac{c_{i}z_{i}}{0.8E_{i,\text{cap}}},\\
\text{s.t., }\gamma_{i}\eta_{i}^{c}[1000\times a_{k}(p_{i,t}^{c})^{2}+n_{i}b_{k}p_{i,t}^{c}]\\
+(1-\gamma_{i})[1000\times a_{k}(p_{i,t}^{d})^{2}+n_{i}b_{k}p_{i,t}^{d}]/\eta_{i}^{d} & \le z_{i},\thinspace\thinspace\forall\, k\in\{1,2,3\},
\end{aligned}
\]
}%
where $c_{i}$ is the unit capital cost (\$/Wh) to purchase the
ESS $i$; $z_i$ is an auxiliary variable; $\gamma_{i}$ is the fraction of a single cyclic aging cost incurred by fully charging the battery from empty; and $n_{i}\triangleq E_{i,\text{cap}}/0.0081$,
which is the number of battery modules that form the ESS $i$, each with a capacity of 0.0081
kW. The battery's price $c_{i}$ is divided by
0.8 because in practice a battery typically needs to be replaced if
its useable capacity drops down to 20\% of its original capacity.
The inequality constraints are essential to enforce the piece-wise linear aging
models, and the coefficient $(a_{k},b_{k})$ takes the
values of (0.0020 h, 0.0086 Wh), (0.0026 h, 0.0060 Wh) and (0.0134
h, -0.0884 Wh) for $k$ equal to 1, 2 and 3, respectively \cite{trippe2014charging}. In the
simulations, we set $c_{1}=c_{2}=0.25$ \$/Wh (which is the ESS price pursued in the next few years \cite{nykvist2015rapidly}),
and $\gamma_{1}=\gamma_{2}=0.5$, and the initial SOCs as $s_{1,0}=0.5$
and $s_{2,0}=0.6$ for the two ESSs with the smaller and the larger
capacities, respectively.

\begin{table}
\protect\caption{\textsc{\small{}Parameters of ESSs}.\label{tb: ES parameters}}

\centering{}{\footnotesize{}}%
\begin{tabular}{cccccccc}
\hline
{\footnotesize{}Type} & {\footnotesize{}$E_{i,\text{cap}}$} & {\footnotesize{}$s_{i,\min}$} & {\footnotesize{}$s_{i,\max}$} & {\footnotesize{}$p_{i,\max}^{c}$} & {\footnotesize{}$p_{i,\max}^{d}$} & {\footnotesize{}$\eta_{i}^{c}$ } & {\footnotesize{}$\eta_{i}^{d}$}\tabularnewline
\hline
{\footnotesize{}1} & {\footnotesize{}480} & {\footnotesize{}0.2} & {\footnotesize{}0.9} & {\footnotesize{}34} & {\footnotesize{}25} & {\footnotesize{}0.82} & {\footnotesize{}0.88}\tabularnewline
{\footnotesize{}2} & {\footnotesize{}720} & {\footnotesize{}0.2} & {\footnotesize{}0.9} & {\footnotesize{}49} & {\footnotesize{}37} & {\footnotesize{}0.85} & {\footnotesize{}0.90}\tabularnewline
\hline
\end{tabular}
\end{table}

\textbf{Energy market: }The price of buying electricity from the energy
market is obtained from PG\&E \cite{energyprice,lu2013online}.
The price has three different values in each day, which is
0.232 \$/kWh during the peak hours 12-18, 0.103 \$/kWh during the
mid-peak hours 8-12 and 18-20, and 0.056 \$/kWh during the remaining off-peak
hours. The price of selling electricity to the energy market is assumed
to be 60\% of the purchase price. The purchase and sale power are both upper
bounded by 1 MW.

\textbf{Load demands: }The demand trace is the hourly electricity consumption data of a college in California
for a typical summer week  \cite{consumptiontrace2002}. Its maximum power is scaled to be 3 MW, as shown in Fig. \ref{fig: One-stage-vs-two-stage}. The inelastic demands are randomly generated within 70\%-90\% of the total demands, and the rest are treated as the elastic demands. The demands are assumed
to be forecasted with zero-mean errors which follow uniform distributions.
The maximum forecast error is proportional
to the absolute change in the corresponding demand of the time slot
w.r.t. that of the previous time slot, and the proportion coefficient
is equal to 5\% and 10\% for the inelastic and elastic demands,
respectively. The error is capped within 80\% of the minimum and 120\% of the maximum demand.

\textbf{Service requirements:} The long-term average percentage of
unsatisfied elastic demand is set to be no more than 30\%, while the
maximum percentage is set to 40\%. When there is a supply shortage
(or surplus), the shortage (or surplus) cost at a price of 0.06 (or
0.07) \$/kWh is imposed to the microgrid as a penalty. In addition,
the total amount of carbon emission is kept below 1,337.6 kg/h (which
is the amount generated by the three CGs working at 95\% of
their maximum output powers), and the operating reserve is maintained
no less than 5\% of the aggregate maximum output power of the three
CGs.

\subsection{Case Studies}

\subsubsection{Value of the second-stage optimization}

\begin{figure}
\noindent \begin{centering}
\includegraphics[scale=0.515]{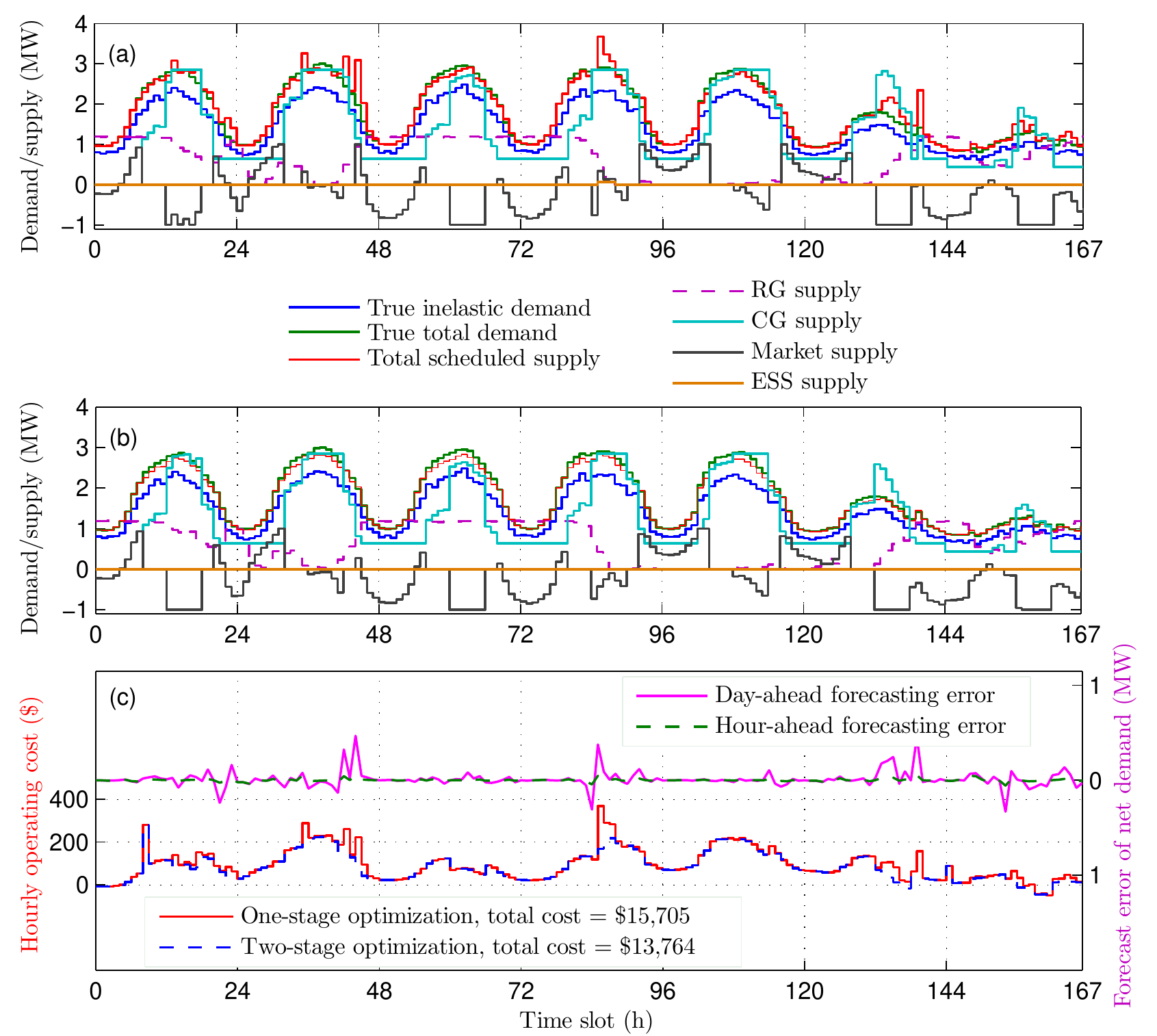}
\par\end{centering}

\begin{centering}
\protect\caption{\small One-stage optimization vs. two-stage optimization: (a) only the first stage of the proposed optimization is applied; (b) both stages of the proposed optimization are applied; and (c) the hourly operating costs and hourly demand forecasting errors. \label{fig: One-stage-vs-two-stage}}
\par\end{centering}
\end{figure}

The real-time economic dispatch
performed in the second stage aims to optimize the amount of energy generated, stored and traded
based on the hour-ahead forecast data, which respects the UC solution obtained in the first stage. The added value comes from the more
accurate data available. This is verified by the simulation
results shown in Fig. \ref{fig: One-stage-vs-two-stage} (the parameters
$\alpha_{\text{avg}}$ and $\alpha_{\max}$ are set to 0.3 to enable a fair comparison). The second-stage optimization
contributes to a lower operating cost compared to that when only
the first-stage optimization is applied. And the cost (\$13,764) is close to the ideal benchmark cost (\$13,537) obtained by optimizing the one-week schedule as a whole using the error-free forecast data.

\subsubsection{Value of having a good model}

\begin{figure*}
\vspace{-12pt}
\noindent \begin{centering}
\includegraphics[width=14cm, height=7.8cm]{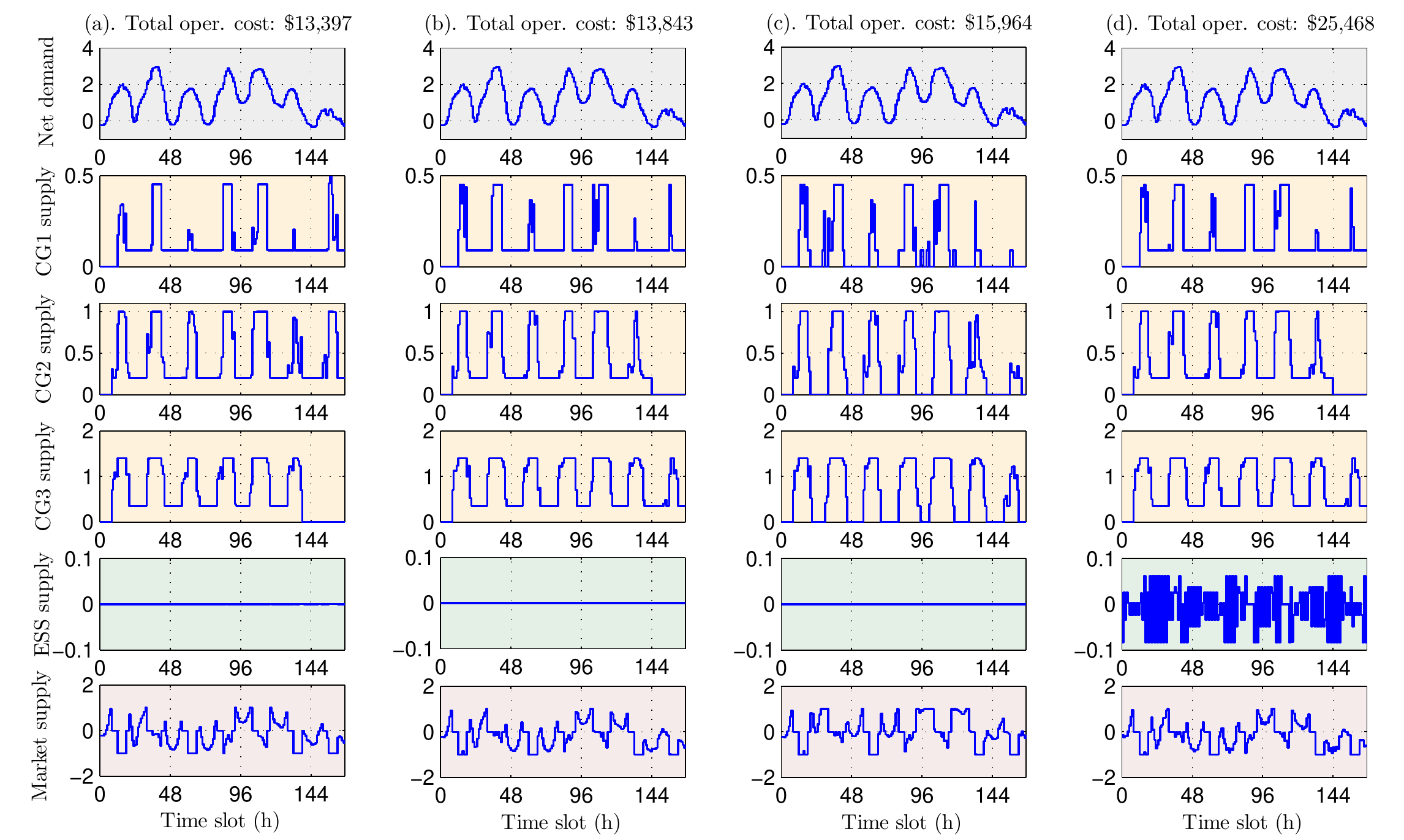}
\par\end{centering}

\begin{centering}
\protect\caption{\small Scheduled generation (MW) based on different optimization models: (a) the ideal benchmark generation; (b) the proposed generation based on an accurate model; (c) the proposed generation without incorporating the CGs' start-up \& shut-down cost model; and (d) the proposed generation without incorporating the ESSs' charging \& discharging cost model.\label{fig: value-of-having-a-good-model}}
\end{centering}
\vspace{-12pt}
\end{figure*}

The start-up/shut-down and charging/discharging costs are often inappropriately modeled to enable simpler models in the literature,
e.g., \cite{sundistributed2015,shireal2015}. It is of interest to examine the consequences. As shown in Fig. \ref{fig: value-of-having-a-good-model},
when the start-up/shut-down cost is not incorporated into the optimization model,
the schedule leads to frequent shut-downs and start-ups of the CGs. This increases
the actual operating cost to \$15,964, which is higher than its counterpart (\$13,843) when both types of costs are incorporated. The consequence is more serious when the charging and discharging costs are not included in the optimization, which leads to frequent charges and discharges of the ESSs and consequently a higher operating cost of \$25,468. The results demonstrate the importance of having appropriate cost models for both CGs and ESSs, which can considerably reduce the total operating cost.

\subsubsection{Sensitivity to forecast errors}

\begin{figure}
\noindent \begin{centering}
\includegraphics[scale=0.53]{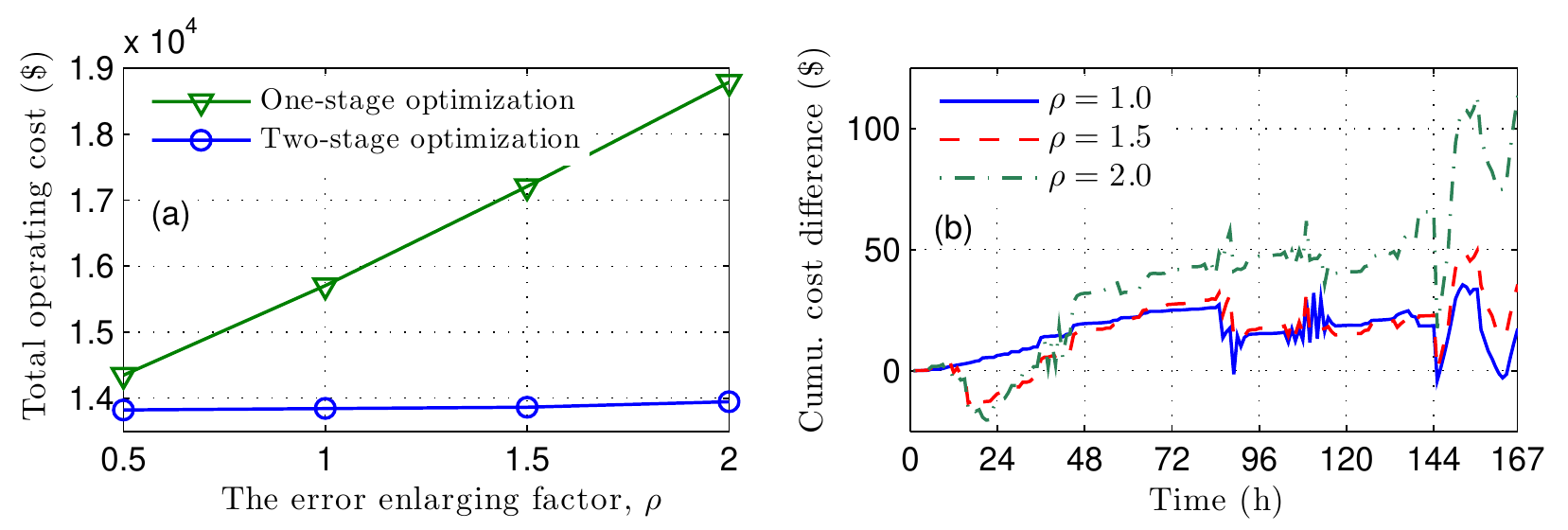}
\par\end{centering}

\centering{}\protect\caption{\small Impact of the forecasting errors: (a) the total operating cost w.r.t. the factor ($\rho$) multiplied to the forecast errors; and (b) the cumulative cost differences for $\rho=1, 1.5$ and 2 relative to $\rho=0.5$. \label{fig: impact-of-forecast-error}}
\vspace{-12pt}
\end{figure}

As both day-ahead and hour-ahead forecast errors of the demands and renewable generation are increased by multiplying a common factor $\rho$, the impacts to the generation schedule are reflected in the increase of the total operating cost. As shown in Fig. \ref{fig: impact-of-forecast-error}(a), the cost increment by applying the two-stage optimization approach is small for $\rho$ increasing from 0.5 to 2. In contrast, the cost increase is significantly much larger if only the first-stage optimization is applied. This indicates that the two-stage approach is robust against forecast errors, mostly owing to the second-stage optimization with a shorter forecast lead time. The curves in Fig. \ref{fig: impact-of-forecast-error}(b) depict the evolutions of the cumulative differences of the total operating costs for $\rho=1, 1.5$ and 2 relative to that for $\rho=0.5$, which verify again that the operating cost increases with the forecast errors.

\section{Conclusions}

This work presented a two-stage optimization approach to schedule
the energy generation in a microgrid based on a realistic model. The approach integrates a day-ahead UC and an hour-ahead economic dispatch as a way to handle large uncertainties in the demand and the renewable generation. The combined solution approximately optimizes the long-term average operating cost of the microgrid subject to operational constraints and both short-term and long-term service requirements. The simulations conducted with real trace data verify the efficacy of the proposed approach, and demonstrate the importance and
value of adopting realistic models for the components, especially the ESSs, in a microgrid.

Our future work includes conducting performance analysis of the proposed scheduling approach and extending the case studies to other aspects of the microgrid, such as the impacts of the service requirements and design
parameters on the generation scheduling and economic performance.

\section*{Acknowledgment}
This work was supported in part by the Energy Innovation Programme Office (EIPO) through the NRF and Singapore EDB. The authors would like to thank Dr. Feng Guo from NEC Laboratories America, Inc. for his useful comments.

\bibliographystyle{IEEEtran}
\bibliography{HuWangGooi_PSCC2016}

\end{document}